\documentclass[letterpaper, 11pt]{amsart}

\usepackage{amsmath,amsthm,amsfonts,amssymb,amscd}
\usepackage{bbm}
\usepackage{bm}
\usepackage{tikz}
\usepackage{tikz-cd}
\usepackage{appendix}
\usepackage{BOONDOX-calo}
\usepackage{standalone}
\usetikzlibrary{arrows,chains,matrix,positioning,scopes, cd}
\usepackage[letterpaper, left=3cm,right=3cm, top=3cm, bottom=3cm]{geometry}
\usepackage{adjustbox}
\usepackage{enumitem}
\usepackage[all]{xy}
\usepackage[
colorlinks=true, citecolor=blue, linkcolor=blue, urlcolor=red]{hyperref}
\usepackage{multicol}
\allowdisplaybreaks
\usepackage{mathrsfs}
\usepackage{here}
\usepackage{blkarray}
\usepackage{caption} 
\usepackage{footnote}
\usetikzlibrary{patterns}
\usepackage{accents}

\newcommand{\R}{{\mathbb R}}

\newcommand{\Z}{{\mathbb Z}}
\newcommand{\Q}{{\mathbb Q}}

\newcommand\dual{\raise0.9ex\hbox{$\scriptscriptstyle\vee$}}

\newcommand{\mathsym}[1]{{}}
\newcommand{\unicode}[1]{{}}

\newcommand{\dd}{\textup{d}}

\newcommand{\D}{{\Delta}}

\theoremstyle{plain}
\newtheorem{thm}{Theorem}

\numberwithin{thm}{section}
\numberwithin{equation}{section}

\newtheorem{manualtheoreminner}{Theorem}
\newenvironment{thm'}[1]{%
  \renewcommand\themanualtheoreminner{#1}%
  \manualtheoreminner
}{\endmanualtheoreminner}

\theoremstyle{definition}

\theoremstyle{remark}

\numberwithin{equation}{section}

\tikzset{>=stealth}

\makeatletter
\def\@seccntformat#1{%
  \protect\textup{\protect\@secnumfont
    \ifnum\pdfstrcmp{subsection}{#1}=0 \bfseries\fi
    \csname the#1\endcsname
    \protect\@secnumpunct
  }%
}  
\makeatother

\makeatletter
\@namedef{subjclassname@2020}{%
  $2020$ Mathematics Subject Classification}
\makeatother

\begin{document}

\title[Improved rational approximations to Catalan's constant]{Improved rational approximations to Catalan's constant}
\author{Payman Eskandari}
\address{Department of Mathematics and Statistics, University of Winnipeg, Winnipeg MB, Canada }
\email{p.eskandari@uwinnipeg.ca}
\begin{abstract}
The works of Rivoal, Zudilin, and Krattenthaler (\cite{RZ},\cite{Zu2}, \cite{Ri1},\cite{KR}) and Nesterenko \cite{Ne} give an explicit sequence $p_n/q_n$ ($p_n,q_n\in \Z, q_n>0$) of rational approximations to Catalan's constant $G=1-1/3^2+1/5^2-1/7^2+1/9^2-\cdots$ satisfying  $\vert G-p_n/q_n\vert \leq 1/q_n^{0.52}$ for sufficiently large $n$. We build on their works to give an explicit sequence $p_n/q_n$ ($p_n,q_n\in \Z, q_n>0$) with $\vert G-p_n/q_n\vert \leq 1/{q}_n^{0.62}$ for sufficiently large $n$.
\end{abstract}
\maketitle
\vspace{-.3in}
\section{Introduction}
\subsection{Context and the statement of the result}
Given a real number $\alpha$ that we hope to prove to be irrational and a sequence of rational approximations
\[\alpha\neq p_n/q_n \rightarrow \alpha, \quad\quad p_n,q_n\in\Z, q_n>0, \gcd(p_n,q_n)=1, n\geq 1 \]
to $\alpha$, it is natural to consider the invariant of the sequence defined by
\begin{equation}\label{eq1}
\begin{split}
\kappa&:=\sup  \left\{\kappa'\in \R: \ \vert \alpha-p_n/q_n\vert \leq 1/q_n^{\kappa'} \ \ \text{for all sufficiently large $n$} \right\}\\
& \ = \liminf_{n\rightarrow \infty} \ \frac{-\log|\alpha-p_n/q_n|}{\log q_n}.
\end{split}
\end{equation}
We will call this the (denominator) \emph{exponent index} of the sequence $p_n/q_n$. If there exists a sequence of rational approximations to $\alpha$ with exponent index $>1$, then $\alpha$ is irrational. It is a classical approach towards irrationality proofs to try to construct a sequence of rational approximations with $\kappa>1$; for example, this is the approach that leads to proofs of irrationality of numbers such as $\alpha=\pi,e$ and perhaps most famously in the last few decades, $\zeta(3)$. In these proofs, the approximating sequence is explicitly defined using continued fraction convergents, partial sums in an infinite sum converging to $\alpha$, linear recurrences with polynomial coefficients (as in Ap\'ery's original proof of irrationality of $\zeta(3)$ \cite{Apery}), or integral expressions (as in Beukers' proof of  irrationality of $\zeta(3)$ \cite{Beukers}). Following Nesterenko \cite{Ne}, let us call constructions of these types \emph{efficient}.

In cases that one cannot prove irrationality of $\alpha$, one can still use \eqref{eq1} to compare the qualities of different efficient rational approximations of $\alpha$, with a sequence of rational approximations of $\alpha$ with a larger exponent index being ``better" than a sequence of approximations of $\alpha$ with a smaller exponent index.\footnote{Note that unlike the notion of convergence rate, the exponent index incorporates both the Archimedean and non-Archimedean properties of the approximating sequence.} The reader can find a summary of the known results in this direction for several classical constants in \cite{Ne}.

This paper is about efficient rational approximations to Catalan's constant
\[
G:=1-\frac{1}{3^2}+\frac{1}{5^2}-\frac{1}{7^2}+\frac{1}{9^2}-\cdots=L(2,\chi_{-4}),
\]
where $L(s,\chi_{-4})$ is the Dirichlet $L$-function for the unique nontrivial character mod 4. As far as the author knows at the time of writing this article, the best (i.e., with largest exponent index) currently available published unconditional efficient rational approximations to Catalan's constant are given by a sequence that first appeared in the work \cite{RZ} of Rivoal and Zudilin. The sequence terms are $v_n/u_n$, where $v_n,u_n$ are well-defined \emph{rational} numbers such that 
\begin{equation}\label{eq2}
   u_nG-v_n = \int\limits_{[0,1]^2} \frac{x^{n-1/2}y^n(1-x)^n(1-y)^{n-1/2}}{(1-xy)^{n+1}}\, \dd x\dd y.
\end{equation}
That such $u_n,v_n\in \Q$ exist is proved in \cite{RZ} (see Lemma 2 and \S 9 therein). The sequence was originally defined in terms of a 
well-poised hypergeometric ${}_{6}F_{5}$ series at $-1$ for some specially chosen parameters depending on $n$, which then by Bailey's transformation can be related to the value of ${}_{3}F_{2}$ at 1 for some suitable parameters, which in turn becomes the integral in \eqref{eq2}. Zudilin \cite{Zu1} has shown that the sequences $v_n,u_n$ satisfy an Ap\'ery-style linear recurrence of order 2. The denominators of $v_n$ and $u_n$ have been carefully studied, and sharp bounds are proved in \cite{KR}, \cite{Ri1}, and later independently in \cite{Ne} (see the last paragraph of \cite{Ne} for the chronology of the developments). Nesterenko has used these denominator bounds\footnote{In fact, Zudilin's slightly weaker bounds from \cite[Theorem 1]{Zu2} suffice for this and the results of the present paper.} to deduce that the sequence $v_n/u_n$ approximates $G$ with an exponent index $\geq 0.52$ (this is the sequence of \cite[Theorem 1]{Ne}, see also \S 4 therein). Experiments suggest that the exponent index of this sequence is indeed roughly equal to the proven bound.

The goal of this paper is to prove the following result:
\begin{thm}\label{thm}
There exists an efficiently constructed sequence of rational approximations of $G$ with exponent index $\kappa\geq 0.62$. More explicitly, there exists an efficiently constructed sequence of rational numbers $\widetilde{p}_n/\widetilde{q}_n$ ($\widetilde{p}_n,\widetilde{q}_n\in \Z, \widetilde{q}_n>0$) with
\[
0<\vert G-\widetilde{p}_n/\widetilde{q}_n\vert \leq 1/\widetilde{q}_n^{\, 0.62}
\]
for all sufficiently large $n$.
\end{thm}
Our experiments suggest that the actual value of $\kappa$ for our sequence is also roughly the proven bound. Note the qualification ``efficient" (in the sense explained earlier) in the statement of the theorem.\footnote{Of course, for any $\alpha\in \R$ there always exists a sequence of rational approximations with exponent index 1.} Our sequence will be defined using an integral construction, as we now explain.

\subsection{Idea of the construction}
Rivoal-Zudilin \cite{RZ} and Nesterenko \cite{Ne} give a more general 5-parameter family of ${}_3F_{2}$-values leading to linear forms in 1 and $G$: adopting Nesterenko's parametrization, for integers $b_1\geq b_2\geq a_1\geq a_2\geq a_3> 0$ with $b_1+b_2\geq a_1+a_2+a_3$, one has
\begin{equation}\label{eq22}
\int\limits_{[0,1]^2}\frac{x^{a_1-1/2}y^{a_2}(1-x)^{b_1-a_1}(1-y)^{b_2-a_2-1/2}}{(1-xy)^{a_3+1}}\, \dd x\dd y = uG-v
\end{equation}
for some well-defined $u,v\in \Q$. Moreover, Nesterenko proves bounds on the lowest denominators of $u,v$ (which improve on the bounds in \cite{RZ}). The special 1-parameter subfamily \eqref{eq2} is recovered by setting $b_2=b_1=2n$ and $a_1=a_2=a_3=n$. The best known denominator bounds for this subfamily, due to Rivoal and Krattenthaler (\cite{KR} and \cite{Ri1}) and independently Nesterenko \cite{Ne}, are 
\begin{equation}\label{eq4}
4^{2n}u_n\in \Z, \quad 4^{2n}L_{2n}^2v_n\in \Z,
\end{equation}
where $u_n,v_n$ are as in \eqref{eq2} and $L_m$ denotes the least common multiple of $1,2,\ldots,m$. These bounds seem to be strikingly sharp: they give a denominator
\[
4^{2n}L_{2n}^2\,|u_n|
\]
for $v_n/u_n$, and experiments suggest that this is almost exactly the lowest denominator of $v_n/u_n$. This is useful, as wasteful denominator bounds make it difficult to accurately estimate the exponent index. On the other hand, the Nesterenko bounds for other subfamilies of \eqref{eq22} do not seem to be necessarily as sharp; see for example, the family in \S 5 of \cite{Ne}, which might have an exponent index of $\approx 0.55$, but this is harder to establish and remains open, because the lowest denominators of the rational approximations seem to be actually smaller than what the established bounds give.

With this in mind, the idea is to remain in the special Rivoal--Zudilin subfamily \eqref{eq2}, but allow ourselves to consider $\Z$-linear combinations of the integrals in the family. The denominators of these combinations will thus still be controlled by the denominator bounds for the integrals \eqref{eq2}. Set
\begin{equation}\label{eq3}
\psi:=\frac{xy(1-x)(1-y)}{1-xy}, \quad\quad J_n:=\int\limits_{[0,1]^2} \psi^n\frac{\dd x\dd y}{\sqrt{x(1-y)}(1-xy)},
\end{equation}
so that $J_n$ is simply \eqref{eq2}. One may now consider integrals of the form
\[
\int\limits_{[0,1]^2} \psi^n \cdot T(\psi)^m \,\frac{\dd x\dd y}{\sqrt{x(1-y)}\,(1-xy)},
\]
where $T\in \Z[t]$ is an auxiliary polynomial and $m\sim \alpha n$, $\alpha$ to be optimized. Such an integral is a $\Z$-linear combination of the $J_{n+i}$ for $0\leq i\leq m\cdot \deg(T)$. Our construction takes the particular choice of an auxiliary polynomial
\[
T(t)=1-16t.
\]
The introduction of the parameter $\alpha$ (where $m\sim \alpha n$) improves the Archimedean decay of the forms at the cost of increasing their denominators. But crucially, this denominator loss is smaller than the naive one, thanks to a $2$-adic saving arising from the factor $16$. As a result, the choice of the auxiliary polynomial $T(t)=1-16t$ results, for any $\alpha>0$, in an improvement of the linear forms compared to the original Rivoal--Zudilin forms \eqref{eq2} (see Remark (2) in \S \ref{sec: remarks} for more details). Optimizing for $\alpha$, the sequence of Theorem \ref{thm} is obtained for $\alpha\approx 3.32$. Needless to say, Theorem \ref{thm} is still far from allowing us to prove irrationality of Catalan's constant.

Note that the idea of considering linear combinations of known linear forms to prove existence of better quality linear forms already appears in Zudilin's \cite{Zu3}. This idea has been recently further developed by Brown in \cite{Br2}. A by-product of the present paper is perhaps the illustration that one can sometimes first try to improve the quality of the linear forms by considering suitable \emph{explicit} linear combinations, before resorting to determinant-based and geometry of numbers existence methods.

\subsection{Further context and related work}
The organization of the paper somewhat unjustly hides the route that resulted in Theorem \ref{thm}. Brown \cite{Br16} has proposed that algebraic geometry, and in particular, geometric constructions of motives with desired periods could lead to the discovery of good rational approximations. The present work is a sequel to \cite{EMN}, where with Murty and Nemoto we initiated Brown's program for $G$, by constructing a 2-dimensional mixed motive with $G$ as a period. A corollary of the construction was that for every symmetric polynomial $F\in \mathbb{Q}[x^2,y^2]$, the integral of $F\dd x\dd y/(1-x^2-y^2)^{k+1}$ over the simplex $\D$ defined by $x,y,1-x-y\geq 0$ is in $\Q+\Q G$, provided that the differential form is ``integrable" over $\D$ (this is to guarantee convergence, see \cite{EMN}). The boundary of $\D$ and the required symmetry lead one to consider integrals of the form
\begin{equation}\label{eq21}
\int_\D \frac{f^n h^m}{g^{k+1}}\dd x\dd y,
\end{equation}
where
\[
f=x^2y^2,\quad g=1-x^2-y^2,\quad h=\prod\limits_{\delta,\epsilon\in\{\pm 1\}} (1-\delta x-\epsilon y)=g^2-4f
\]
and $k\leq 2m+2n$ (for convergence). Our original construction of the rational approximations of Theorem \ref{thm} was in the framework of the integrals \eqref{eq21}, in terms of suitable linear combinations of the integrals \eqref{eq21} with $m=n$, $k=2n$, rather than linear combinations of the hypergeometric integrals $J_n$. It turns out that one has the perhaps rather surprising identity 
\begin{equation}\label{eq20}
\int\limits_\D \frac{f^nh^n}{g^{2n+1}}\dd x\dd y = 2^{-2n-3} \!\!\!
\int\limits_{[0,1]^2} \frac{x^{n-1/2}(1-x)^ny^n(1-y)^{n-1/2}}{(1-xy)^{n+1}}\dd x\dd y.
\end{equation}
Sharp denominator bounds were already proved in the hypergeometric framework of the right hand side by the works earlier mentioned. We then realized that the entire construction can be done in the hypergeometric framework, which led to the current presentation of this paper. The identity \eqref{eq20} is a special case of a more general relation between the two families \eqref{eq21} and \eqref{eq22}. See \cite{Es2} for more details.

Finally, after completing this work we learned that Carlo Viola and Raffaele Marcovecchio \cite{MaViMFO} have obtained rational approximations to $G$ of a similar quality to the ones in this paper.\footnote{I thank Wadim Zudilin for bringing this to my attention.} The constructions of this paper are different from and somewhat simpler than those of Viola and Marcoveccio.

\section{Constructions of the linear forms}
Throughout the paper, $\psi$ and $J_n$ are as in \eqref{eq3}. We reserve the notation $u_n,v_n$ for the rational numbers satisfying
\[J_n=u_nG-v_n,\]
explicitly given in \cite{Ne}. By Theorems 1 and 3 of \cite{Ne}, the rational numbers $u_n,v_n$ satisfy \eqref{eq4} (also see \S 4 therein). Nesterenko gives the following explicit formula for $u_n$:
\begin{equation}\label{eq23}
    u_n = 8(-1)^{n}4^{-2n} \sum_{j=0}^n 4^j\binom{n}{j} \frac{(2n+2j)!j!}{(n+j)!(2j)!n!} \frac{(2n)!j!}{n!(2j)!(n-j)!}.
\end{equation}
Indeed, our $u_n$ is Nesterenko's $4B(1)$ in \cite{Ne} for the choice of parameters $a_3=a_2=a_1=n$ and $b_1=b_2=2n$. The formula above is what Theorem 2 and equation (3.10) of \cite{Ne} give. (Notes: The choice of parameters $a_k,b_\ell$ forces one to always be in the second case of the computations of the proof of Theorem 3 of \cite{Ne}. More importantly, note that there is a sign typo in the formula (3.10) of \cite{Ne}; the sign typo can be seen by computing the left-hand side of (3.8) of said reference by means of (2.8) therein.)

Focusing on the inner sum in \eqref{eq23} for the moment, the initial term (i.e., with $j=0$) is 
\[\binom{2n}{n}^{\!2}.\] 
The quotient of the $(j+1)$-term by the $j$-term simplifies to
\[
\frac{(j-n)^2(j+n+1/2)}{(j+1/2)^2}\frac{1}{j+1},
\]
which is exactly the corresponding successive quotient of the finite hypergeometric sum
\begin{equation}\label{eq7}
F_n:={}_3F_2\left(\begin{matrix}-n,-n,n+\frac12\\[1mm]
\frac12,\frac12\end{matrix};1\right).
\end{equation}
Thus
\begin{equation}\label{eq6}
u_n = 8(-1)^{n}4^{-2n}\binom{2n}{n}^{\!2}  F_n.
\end{equation}

We now define and set some notation for the linear forms leading to Theorem \ref{thm}. 
For $m,n\geq 0$, we shall set
\begin{equation}\label{eq30}
\begin{split}
    \widetilde{J}_{n,m} &:= \int\limits_{[0,1]^2} \psi^n \cdot (1-16\psi)^m \,\frac{\dd x\dd y}{\sqrt{x(1-y)}\,(1-xy)}\\
    &= \sum_{i=0}^m (-16)^i{m\choose i} \, J_{n+i}.
\end{split}
\end{equation}
We will treat $m$ as a function of $n$ with $m/n \rightarrow \alpha $ as $n\rightarrow\infty$, where $\alpha\in [0,\infty)$ is a fixed number which will be chosen later. We will write
\[
\widetilde{J}_{n,m} = \widetilde{u}_{n,m}G - \widetilde{v}_{n,m},
\]
where
\begin{equation}\label{eq13}
\widetilde{u}_{n,m} := \sum_{i=0}^m (-16)^i\binom{m}{i} \, u_{n+i}, \quad\quad 
\widetilde{v}_{n,m} = \sum_{i=0}^m (-16)^i\binom{m}{i} \, v_{n+i}.
\end{equation}
Rewriting $u_{n+i}$ using \eqref{eq6}, after simplification we get
\begin{equation}\label{eq8}
\widetilde{u}_{n,m} = 8(-1)^{n}4^{-2n} \sum_{i=0}^m \binom{m}{i} \binom{2n+2i}{n+i}^{\!2} F_{n+i}.
\end{equation}
Notice the role played here by the factor $(-16)^i$ of the summands of $\widetilde{u}_{n,m}$ in \eqref{eq13}. We shall set
\begin{equation}\label{eq11}
\widetilde{F}_{n,m}:=\sum_{i=0}^m \binom{m}{i} \binom{2n+2i}{n+i}^{\!2} F_{n+i}
\end{equation}
so that
\begin{equation}\label{eq24}
\widetilde{u}_{n,m} = 8(-1)^{n}4^{-2n} \widetilde{F}_{n,m}.
\end{equation}
Note that since each $F_{n}$ is positive (as it is clear from the definition \eqref{eq7}), $\widetilde{F}_{n,m}$ is also positive. Thus $\widetilde{u}_{n,m}$ is nonzero.

Turning our focus to the denominators of $\widetilde{u}_{n,m}$ and $\widetilde{v}_{n,m}$, by the bounds \eqref{eq4}, each of $4^{2n+2i}u_{n+i}$ and $4^{2n+2i}L^2_{2n+2i}v_{n+i}$ are integers. Thus by \eqref{eq13},
\begin{equation}\label{eq5}
4^{2n}\widetilde{u}_{n,m}\in \Z, \quad\quad 4^{2n}L^2_{2n+2m}\widetilde{v}_{n,m}\in\Z.
\end{equation}
Again notice the role played here by the factor $16^i$ of the summands of \eqref{eq13}.

\section{Asymptotics of $\widetilde{u}_{n,m}$}\label{sec: asymptotics of utilde}
Here we will study the asymptotic behaviour of $\widetilde{u}_{n,m}$ to the extent needed for the purposes of this paper. Consider the terminating hypergeometric value $F_n$ defined in \eqref{eq7}. Write
\[
F_n=\sum_{j=0}^n a_{n,j},\qquad
a_{n,j}=\frac{(-n)_j^2(n+\frac12)_j}{(\frac12)_j^2 \,j!}=
\frac{\Gamma(n+1)^2\Gamma(n+j+\frac12)\Gamma(\frac12)^2}
{\Gamma(n-j+1)^2\Gamma(n+\frac12)\Gamma(j+\frac12)^2\Gamma(j+1)},
\]
where $(a)_j=\Gamma(a+j)/\Gamma(a)$ is the Pochhammer symbol. Since every $a_{n,j}$ (with $0\leq j\leq n$) is positive, we can try to derive the asymptotics of $F_n$ using the discrete analogue of the Laplace method. For the purpose of this paper, we shall only need the asymptotics of $(\log F_n)/n$, which can be handled in a slightly more elementary way, as follows. Following Nesterenko \cite{Ne}, we start with following version of Stirling's formula (due to Binet):
\begin{equation}\label{eq Sterling}
\left|\log\Gamma(x) - x\log x +x +\frac{1}{2}\log x-\frac{1}{2}\log(2\pi)\right| \leq K/x,
\end{equation}
valid for all $x>0$, with $K>0$ being an absolute constant (\cite{WW}, page 249). For any $0\leq j\leq n$, we use this to estimate each term in
\[
\begin{split}
    \log  a_{n,j} = &2\log\Gamma(n+1)+\log\Gamma(n+j+\frac{1}{2})-2\log\Gamma(n-j+1)\\
    &-\log\Gamma(n+\frac{1}{2})-2\log\Gamma(j+\frac{1}{2})-\log\Gamma(j+1)+ O(1),
\end{split}
\]
where the constant in $O(1)$ is absolute. We obtain
\[\begin{split}
    \log  a_{n,j} &\stackrel{(\ast)}{=} n\log n+(n+j)\log(n+j)-2(n-j)\log(n-j)-3j\log j+ O(\log n)\\
    &\stackrel{(\ast\ast)}{=} (n+j)\log(1+\frac jn)-2(n-j)\log(1-\frac jn)-3j\log (\frac jn)+ O(\log n),
\end{split}
\]
where the constant in $O(\log n)$ is absolute and here and in what follows, $0\log 0$ is to be interpreted as $0$. (For $(\ast)$, the sum of the linear contributions coming from different $\Gamma$-factors is an absolute constant, so it is absorbed in the error term. For $(\ast\ast)$, note that the sum of the coefficients of the expression preceding $(\ast\ast)$ is zero, so we may divide the argument of each logarithmic term by $n$.) We thus have
\begin{equation}\label{eq26}
\frac{1}{n}\log a_{n,j} = (1+\frac jn)\log(1+\frac jn)-2(1-\frac{j}{n})\log(1-\frac jn)-3\frac{j}{n}\log (\frac jn)+ O(\frac{\log n}{n}).
\end{equation}
Consider the function $f$ defined on $(0,1)$ by
\[
f(r) = (1+r)\log(1+r)-2(1-r)\log(1-r)-3r\log r,
\]
and extended continuously to $[0,1]$. The main term of $\frac{1}{n}\log a_{n,j}$ is $f(j/n)$. In view of
\[
f'(r)=\log\!\left(\frac{(1-r)^2(1+r)}{r^3}\right),\qquad f''(r)=\frac{-(r+3)}{r(1-r)(1+r)},
\]
we see that $f$ has a unique critical point in $(0,1)$ at
\[
r_0=\frac{\sqrt 5-1}{2},
\]
and that $f$ attains its unique maximum on $[0,1]$ at $r_0$. Set
\[
\theta:=\frac{1+\sqrt5}{2}=r_0^{-1}.
\]
Using $r_0^2+r_0-1=0$ we get $f(r_0)=5\log\theta$. Since the $j/n$ ($0\leq j\leq n$) can get arbitrarily close to $r_0$ as $n\rightarrow \infty$, from \eqref{eq26} we obtain
\begin{equation}\label{eq25}
    \lim_{n\rightarrow \infty} \max_{0\leq j\leq n} \frac{1}{n}\log a_{n,j} = f(r_0)=5\log\theta.
\end{equation}
Since every $a_{n,j}$ is positive, we have
\[
\max_{0\leq j\leq n}  a_{n,j} \leq F_n\leq (n+1)\max_{0\leq j\leq n}  a_{n,j}.
\]
Thus by \eqref{eq25}, we have
\[
\lim_{n\rightarrow\infty}\frac{1}{n}\log F_n = \lim_{n\rightarrow \infty} \max_{0\leq j\leq n} \frac{1}{n}\log a_{n,j} =5\log\theta.
\]
Combining this with the fact that $\binom{2n}{n}\sim 4^n/\sqrt{\pi n}$, we get
\begin{equation}\label{eq10}
\lim_{n\rightarrow \infty} \left(\binom{2n}{n}^2F_n\right)^{1/n} = R:=16\theta^5=88+40\sqrt5.
\end{equation}

We now use this to estimate $\widetilde{F}_{n,m}$ and hence $\widetilde{u}_{n,m}$. Fix $\varepsilon$ with $0<\varepsilon<R$. There is $N>0$ such that for all $n'>N$,
\[
(R-\varepsilon)^{n'}<\binom{2n'}{n'}^2F_{n'} <(R+\varepsilon)^{n'}.
\]
Taking $n>N$ and $m\geq 0$, multiplying the above with $n'=n+i$ by $\binom{m}{i}$ and summing over $0\leq i\leq m$, on recalling the definition of $\widetilde{F}_{n,m}$ from \eqref{eq11} we obtain
\[
(R-\varepsilon)^n\sum\limits_{i=0}^m \binom{m}{i} (R-\varepsilon)^{i}<\widetilde{F}_{n,m}< (R+\varepsilon)^n\sum\limits_{i=0}^m \binom{m}{i} (R+\varepsilon)^{i},
\]
or equivalently,
\[
(R-\varepsilon)^n (\!1+R-\varepsilon)^m
<\widetilde{F}_{n,m}
< (R+\varepsilon)^n (1+R+\varepsilon)^m.
\]
Considering $m$ as a function of $n$ with $m/n\rightarrow \alpha$ as $n\rightarrow\infty$, we get
\[
    \log(R-\varepsilon) + \alpha\log(1+R-\varepsilon) \leq \liminf_{n\rightarrow\infty} \frac{\log \widetilde{F}_{n,m}}{n}\leq
    \limsup_{n\rightarrow\infty} \frac{\log \widetilde{F}_{n,m}}{n} \leq \log(R+\varepsilon) + \alpha\log(1+R+\varepsilon)
\]
for all $0<\varepsilon<R$. Letting $\varepsilon\rightarrow 0$, we obtain
\[
\lim_{n\rightarrow\infty }\frac{1}{n}\log \widetilde{F}_{n,m} =\log R + \alpha\log(1+R).
\]
Finally, on recalling \eqref{eq24}, we have
\begin{equation}\label{eq12}
\lim_{n\rightarrow\infty }\frac{1}{n}\log|\widetilde{u}_{n,m}| = -\log16 +\log R + \alpha\log(1+R),
\end{equation}
where $m/n\rightarrow \alpha$ as $n\rightarrow \infty$ (and $R$ is as in \eqref{eq10}).

\section{Upper bound for $\widetilde{J}_{n,m}$}\label{sec: asymptotics of Jtilde}
We will denote the sup norm on $[0,1]^2$ by $\|~\|$. For now, fix integers $n>0$ and $m\geq 0$. From the original definition of $\widetilde{J}_{n,m}$ in \eqref{eq30}, we have
\begin{equation}\label{eq29}
    |\widetilde{J}_{n,m}| \leq \|\psi^n(1-16\psi)^m\| \int\limits_{[0,1]^2}\frac{\dd x\dd y}{\sqrt{x(1-y)}\,(1-xy)}= J_0 \cdot \|\psi^n(1-16\psi)^m\|,
\end{equation}
where $J_0=8G<8$ (see \cite{RZ}, \S 9). The reader can check that
\[
\|\psi\|= \max_{0\leq x,y\leq 1} \frac{xy(1-x)(1-y)}{(1-xy)} =\left(\frac{\sqrt{5}-1}{2}\right)^5=\theta^{-5}\approx 0.09017.
\]
The values of $\psi^n(1-16\psi)^m$ on $[0,1]^2$ are the same as the values of $g(t)=t^n(1-16t)^m$ on $[0,\|\psi\|]$. The critical points of $g$ in $[0,\|\psi\|]$ are among $0, n/(16(m+n)), 1/16$ (note that $1/16<\|\psi\|$). The maximum of $|t^n(1-16t)^m|$ on $[0,\|\psi\|]$ is attained at either $n/(16(m+n))$ or the end point $\|\psi\|$. Thus for $m>0$,
\begin{equation}\label{eq27}
\|\psi^n(1-16\psi)^m\| = \max\left\{16^{-n}\left(\frac{n}{m+n}\right)^{\!n}\left(\frac{m}{m+n}\right)^{\!m}
\, , \,\|\psi\|^n (16\|\psi\|-1)^m\right\}
\end{equation}
(note that $16\|\psi\|-1\approx 0.44272>0$). Setting $0^0:=1$, the formula is also valid for $m=0$. 

Now as before, we treat $m$ as a function of $n$ with $m/n\rightarrow \alpha$ as $n\rightarrow\infty$. First let $\alpha>0$. Taking $n$-th roots of \eqref{eq27} and letting $n\rightarrow\infty$, we get
\begin{equation}\label{eq28}
\|\psi^n(1-16\psi)^m\|^{1/n} \rightarrow \max\left\{\frac{1}{16(1+\alpha)}\left(\frac{\alpha}{1+\alpha}\right)^{\!\alpha} \, , \, \|\psi\|\, (16\|\psi\|-1)^\alpha \right\}.
\end{equation}
Defining $(\alpha/(1+\alpha))^\alpha$ at zero by requiring right-continuity at zero (i.e., with value 1), the formula \eqref{eq28} is also valid for $\alpha=0$.

We shall set
\[
M(\alpha) := \max\left\{\frac{1}{16(1+\alpha)}\left(\frac{\alpha}{1+\alpha}\right)^{\!\alpha} \, , \, \|\psi\| \,(16\|\psi\|-1)^\alpha \right\}
\]
on $[0,\infty)$. We will refer to the first and second expressions on the right respectively as the interior and end point terms (referring to what they represent in the variable $t$). They are not difficult to compare as functions of $\alpha$ using standard calculus techniques. Indeed, denoting the interior (resp. end point) term tentatively by $I(\alpha)$ (resp. $E(\alpha)$), we have
\[
\frac{\dd}{\dd\alpha}\log(I/E) = \log(\frac{\alpha}{1+\alpha})-\log(16\|\psi\|-1),
\]
which is negative on $(0,c)$ and positive on $(c,\infty)$ for some $c\in (0,\infty)$ (recall that $16\|\psi\|-1\approx 0.44272<1$). Thus $\log(I/E)$ and hence $I/E$ decreases on $[0,c]$ and increases on $[c,\infty)$. Near zero, $I/E<1$ and for sufficiently large $\alpha$, we have $I/E>1$, so that $E=I$ at exactly one point $\alpha_0\in(0,\infty)$. Numerically, $\alpha_0= 3.31765...$. The end point term dominates for $\alpha$ between $0$ and $\alpha_0= 3.31765...$, at which point the interior point term takes over:
\begin{equation}\label{eq19}
M(\alpha) = \begin{cases}
\|\psi\| \,(16\|\psi\|-1)^\alpha\approx 0.09017\cdot 0.44272^\alpha  \quad\quad &(\alpha\leq\alpha_0= 3.31765...)\\
\frac{1}{16(1+\alpha)}\left(\frac{\alpha}{1+\alpha}\right)^{\!\alpha}\sim \frac{1}{16e\alpha} &(\alpha\geq\alpha_0).
\end{cases}
\end{equation}
Both branches of the function are strictly decreasing on $[0,\infty)$.

By \eqref{eq29} and \eqref{eq28}, we have
\begin{equation}\label{eq15}
\vert\widetilde{J}_{n,m}\vert^{1/n} < M(\alpha)+o(1)
\end{equation}
as $n\rightarrow \infty$ and $m/n\rightarrow \alpha\in [0,\infty)$.

\section{Proof of Theorem \ref{thm}}
We can now establish Theorem \ref{thm}. The sequence of the theorem is the sequence $\widetilde{v}_{n,m}/\widetilde{u}_{n,m}$, where $m$ is a function of $n$ with $m/n\rightarrow \alpha$ as $n\rightarrow \infty$, with $\alpha\in [0,\infty)$ to be chosen. For any choice of $\alpha$, it is possible to make sure that $m$ is always even (simply take $m=\lfloor \alpha n\rfloor$ or $m=\lfloor \alpha n\rfloor+1$ depending on which one is even). By \eqref{eq30}, this will guarantee that $\widetilde{J}_{n,m}>0$. We have
\[
0\neq G-\frac{\widetilde{v}_{n,m}}{\widetilde{u}_{n,m}}=\frac{\widetilde{J}_{n,m}}{\widetilde{u}_{n,m}}\rightarrow 0
\]
as $n\rightarrow \infty$ (and $m/n\rightarrow \alpha$) by our bounds in the last two sections.

Recall from \eqref{eq12} that
\begin{equation}\label{eq14}
\frac{1}{n}\log|\widetilde{u}_{n,m}| = \log R+\alpha\log(1+R)-\log 16 + o(1)
\end{equation}
as $n\rightarrow\infty$. The constant $R$ is given in \eqref{eq10}.

On the other hand, by \eqref{eq15},
\begin{equation}\label{eq17}
\frac{1}{n}\log \widetilde{J}_{n,m} < \log M(\alpha) + o(1)
\end{equation}
as $n\rightarrow \infty$. The function $M(\alpha)$ is as in \eqref{eq19}.

Let us write $\widetilde{v}_{n,m}/\widetilde{u}_{n,m}$ as $\widetilde{p}_{n,m}/\widetilde{q}_{n,m}$ in reduced form (hence in particular, $\widetilde{p}_{n,m}$ and $\widetilde{q}_{n,m}$ are integers). By \eqref{eq5}, 
\[
\widetilde{q}_{n,m} \leq 4^{2n}L^2_{2n+2m} \vert\widetilde{u}_{n,m}\vert \in \Z. 
\]
Combining this with the prime number theorem $\log L_n\sim n$ and \eqref{eq14}, we get that as $n\rightarrow\infty$,
\begin{equation}\label{eq16}
\frac{1}{n} \log \widetilde{q}_{n,m} \leq 4(1+\alpha) \, + \log R+\alpha\log(1+R)+o(1).
\end{equation}
Combining \eqref{eq14}-\eqref{eq16} we get
\[
\begin{split}
\frac{-\log|G-\widetilde{p}_{n,m}/\widetilde{q}_{n,m}|}{\log \widetilde{q}_{n,m}} &= \frac{\log \vert\widetilde{u}_{n,m}\vert-\log \widetilde{J}_{n,m}}{\log \widetilde{q}_{n,m}}\\
&\geq \frac{\log R+\alpha\log(1+R)-\log 16-\log M(\alpha)+o(1)}{4(1+\alpha) \, + \log R+\alpha\log(1+R)+o(1)}\\
&\rightarrow \frac{\log R+\alpha\log(1+R)-\log 16-\log M(\alpha)}{4(1+\alpha) \, + \log R+\alpha\log(1+R)}
\end{split}
\]
as $n\rightarrow\infty$. Thus the exponent index $\kappa$ of our sequence satisfies
\begin{equation}\label{eq18}
\kappa \geq \frac{\log R+\alpha\log(1+R)-\log 16-\log M(\alpha)}{4(1+\alpha) \, + \log R+\alpha\log(1+R)}.
\end{equation}

The final step is now to optimize for $\alpha$. This can be done using standard calculus techniques. Leaving the details out, we see that the function on the right side is increasing on $[0,\alpha_0]$, where $\alpha_0$ is the unique crossing point of the two branches of $M(\alpha)$ (see \eqref{eq19}), and it is decreasing on $[\alpha_0,\infty)$. Thus the right hand side of \eqref{eq18} admits its maximum on $[0,\infty)$ at 
\[
\alpha_0 = 3.31765...,
\]
at which point \eqref{eq18} reads
\[
\kappa \geq 0.62335...
\]

\section{Remarks}\label{sec: remarks}
We end the paper with three remarks.
\medskip\par 
(1) Denoting the right hand side of \eqref{eq18} by $\phi(\alpha)$, for comparison:
\[
\phi(0) \approx 0.52427 ,\ \  \phi(1)  \approx 0.58875 ,\ \  \phi(3)  \approx 0.62097,\ \  \phi(4)  \approx 0.618910, \ \  \lim_{\alpha\rightarrow\infty }\phi(\alpha)  \approx 0.56447.
\]
The case $\alpha=0$ simply recovers the estimate in \cite{Ne} for the Rivoal--Zudilin sequence (note that $\widetilde{J}_{n,0}=J_n$). The parameter $\alpha$ affects both asymptotics of the error term $G-\widetilde{p}_{n,m}/\widetilde{q}_{n,m}$ and the denominator $\widetilde{q}_{n,m}$. Looking at the numerator and denominator of $\phi$ in \eqref{eq18}, increasing $\alpha$ will increase both the numerator (as $0<M(\alpha)<1$ and $M(\alpha)$ is strictly decreasing) and the denominator of $\phi$. What we gain from the improved asymptotics of $G-\widetilde{p}_{n,m}/\widetilde{q}_{n,m}$ dominates the loss in the denominator $\widetilde{q}_{n,m}$ as $\alpha$ goes from $0$ to $\alpha_0$ and hence $\kappa$ (or at least, our estimated value of it) improves. After $\alpha=\alpha_0$, the pattern reverses.
\medskip\par 
(2) Although the denominator $\widetilde{q}_{n,m}$ increases as $\alpha$ increases (with fixed $n$), the 2-adic saving in the denominators $4^{2n}L_{2n+2m}^2$ and $4^{2n}$ of $\widetilde{v}_{n,m}$ and $\widetilde{u}_{n,m}$ compared to the bounds $4^{2n+2m}L_{2n+2m}^2$ and $4^{2n+2m}$ for the denominators of an arbitrary $\Z$-linear combination of Rivoal--Zudilin's $u_N$, $v_N$ with $N\leq n+m$ is crucial in our construction. To make this role more explicit, we can replace 16 in the definition of $\widetilde{J}_{n,m}$ by a positive real variable $z$:
\[
    \widetilde{J}_{n,m}(z) := \int\limits_{[0,1]^2} \psi^n \cdot (1-z\psi)^m \,\frac{\dd x\dd y}{\sqrt{x(1-y)}\,(1-xy)}.
\]
We may now define $\widetilde{u}_{n,m}(z)$ by 
\[
\widetilde{u}_{n,m}(z) := \sum_{i=0}^m (-z)^i\binom{m}{i} \, u_{n+i},
\]
i.e., by replacing $16$ by $z$ in \eqref{eq13}. The contents of \S \ref{sec: asymptotics of utilde} and \S \ref{sec: asymptotics of Jtilde} can be adjusted and we obtain bounds for $\frac{1}{n}\log|\widetilde{u}_{n,m}(z)|$ and $\frac{1}{n}\log|\widetilde{J}_{n,m}(z)|$ which vary continuously with respect to $z$ and specialize to \eqref{eq14} and \eqref{eq17} when $z=16$. In particular, \eqref{eq12} will be modified to
\[
\lim_{n\rightarrow\infty }\frac{1}{n}\log|\widetilde{u}_{n,m}(z)| = -\log16+\log R+\alpha\log(1+zR/16)
\]
where $m/n\rightarrow\alpha$ as $n\rightarrow\infty$.
For an arbitrary positive \emph{integer} $z$, an upper bound for $\widetilde{q}_{n,m}(z)$ (interpreted as the lowest denominator of $\widetilde{v}_{n,m}(z)/\widetilde{u}_{n,m}(z)$ since $z\in\Z$) is
\[
4^{2n+2m}2^{-\overline{\nu}_2(z)m}L_{2n+2m}^2\vert\widetilde{u}_{n,m}(z)\vert,
\]
where $\overline{\nu}_2(z)=\min(\nu_2(z),4)$ with $\nu_2$ being the 2-adic valuation. Our denominator bound \eqref{eq16} now becomes
\[
\frac{1}{n}\log \widetilde{q}_{n,m} \leq 
\alpha(1-\frac{\overline{\nu}_2(z)}{4})\log 16 + 4(1+\alpha) \, + \log R+\alpha\log(1+zR/16)+ o(1).
\]
The role of the $2$-adic valuation of the choice $z=16$ is now evident, as at that point the first term on the right hand side vanishes. For example, if we take $z=17$, we get an extra $\alpha\log 16$ in the denominator of the modified right hand side of \eqref{eq18}, while the rest of the terms in \eqref{eq18} only change slightly, as the bounds in \eqref{eq14} and \eqref{eq17} vary continuously in an Archimedean sense with respect to $z$.
\medskip\par 
(3) Experiments suggest that the lowest denominators of $\widetilde{u}_{n,m}$ and $\widetilde{v}_{n,m}$ are almost exactly $4^{2n}$ and $4^{2n}L_{2n+2m}^2$. Similarly, the bound $\phi(\alpha)$ for the exponent index of a sequence $\widetilde{p}_{n,m}/\widetilde{q}_{n,m}$ with $m\sim \alpha n$ is likely to be almost the exact value.

\section*{Acknowledgments}
This work is a continuation of the joint work \cite{EMN} with Kumar Murty and Yusuke Nemoto. I would like to express my utmost gratitude to Murty and Nemoto for very helpful discussions over the course of this work. I would also like to thank Michel Waldschmidt and Daniel Bertrand for helpful discussions during a visit to Paris and their hospitality. I am also grateful to Francis Brown for helpful correspondence, and in particular, for bringing \cite{Ne} to my attention and generously sharing an early draft of \cite{Br2}. Last but not least, I am thankful to my student Mithuna Threz Arul Milton, who wrote a SageMath code that implemented the calculations of \cite{EMN} for the coefficients of the linear forms of \cite{EMN}. Without the experiments run using her code, this work would not have been possible.

\end{document}